%&amstex
\input amstex
\input amsppt.sty
\define\pro{\prod\limits}
\define\su{\sum\limits}
\define\inte{\int\limits}
\topmatter
\title 
Biorthogonal Ensembles
\endtitle
\author Alexei Borodin
\endauthor
\thanks 
Supported by the Russian Program for Support of Scientific Schools (grant
96-15-96060).\endthanks
\abstract
One object of interest in random matrix theory is a family of point ensembles
(random point configurations) related to various systems of classical
orthogonal polynomials. The paper deals with a one--parametric
deformation of these ensembles, which is defined in terms of the
biorthogonal polynomials of Jacobi, Laguerre and Hermite type.

Our main result is a series of explicit expressions for the correlation
functions in the scaling limit (as the number of points goes to
infinity). As in the classical case, the correlation functions have
determinantal form. They are given by certain new kernels which are
described in terms of the Wright's generalized Bessel function and can be
viewed as a generalization of the well--known sine and Bessel
kernels. 

In contrast to the conventional kernels, the new kernels are
non--symmetric. However, they possess other, rather surprising,
symmetry properties. 

Our approach to finding the limit kernel also differs from the
conventional one, because of lack of a simple explicit
Christoffel--Darboux formula for the biorthogonal polynomials.  
\endabstract
\endtopmatter
\document
\head 1. Introduction \endhead
Orthogonal polynomial ensembles are widely known.  They are characterized by the property that the joint probability density of an $N$-point ensemble has  the form
$$
p(x_1,\ldots,x_N)=const\cdot\pro_{i=1}^N\omega(x_i)\pro_{i<j}(x_i-x_j)^2
$$
for some positive weight function $\omega(x)$. The phase space $I$ of an orthogonal ensemble is a  finite or
infinite interval of the real line. Usually the couple $(\omega(x),I)$ corresponds to one of the classical systems of 
orthogonal polynomials.

These ensembles play a very important role in the random matrix theory, see [Me]. They also serve as a rich source
 of various mathematical problems which include Selberg integrals,
  differential equations for Fredholm determinants 
 and many others, see, for example, [Me], [TW]. 

One of the main properties of orthogonal ensembles is the existence of simple
formulas for all correlation functions.
Namely, the $n$th correlation function has the form
$$
\rho_{nN}(x_1,\ldots,x_n)=\pro_{i=1}^n\omega(x_i)\cdot\det\left [K_N(x_i,x_j)\right]_{i,j=1}^n
\tag *
$$
where  $K_N(x,y)$ is the $N$th Christoffel-Darboux kernel;
$$
K_N(x,y)=\frac{p_N(x)p_{N-1}(y)-p_{N-1}(x)p_{N}(y)}{x-y};  
$$
$\{p_k(x)\}$ is the system of orthogonal polynomials on $I$ with weight
 $\omega(x)$ (see [D], [Me], [NW1]).
 Determinantal formulas (*) often allow to study the asymptotic behaviour
  of an $N$-point biorthogonal ensemble
when the number of points $N$ goes to infinity. It turns out that after the appropriate scaling the correlation functions 
tend to a limit. The limit correlation functions also have determinantal
form with a certain limit kernel. 

The well--known {\it sine kernel} 
$$
\frac{\sin\pi(x-y)}{\pi(x-y)}.
$$
 arises in the
scaling limit of the classical polynomial ensembles in the bulk of
spectrum. In particular, the sine kernel arises in the case of the Hermite
weight function, see, e.g., [Me, section 5.2].

At the left edge of the spectrum, both the Jacobi and Laguerre
weights produce the same {\it Bessel kernel} ([F], [NW2])
$$
\frac{\varphi_1(x)\varphi_2(y)-\varphi_1(y)\varphi_2(x)}{x-y}
$$
where
$$
\varphi_1(x)=J_{\alpha}(2\sqrt{x}),\quad \varphi_2(x)=x\varphi_1'(x);
$$
$\alpha>-1$, $J_{\alpha}(z)$ is the Bessel function.
 
Both the sine and the Bessel
kernels are symmetric and represent positive selfadjoint integral
operators in $L^2$ on $\Bbb R_+$ and on
$\Bbb R$ respectively. 

In this paper we study a one parameter generalization of the orthogonal polynomial ensembles,
 the joint probability densities of our ensembles have the form  
$$
p^b(x_1,\ldots,x_N)=const\cdot\pro_{i=1}^N\omega(x_i)\pro_{i<j}\left[(x_i-x_j)(x_i^{\theta}-x_j^{\theta})\right]
$$
where $\theta$ is a fixed positive number. We call these ensembles {\it biorthogonal}.
 Orthogonal polynomial ensembles correspond to $\theta=1$.
  
Biorthogonal ensembles inherit several nice properties from the orthogonal ones. 
For example, their correlation functions also have determinantal form, and corresponding 
kernels can be expressed via so-called biorthogonal polynomials, see [K1] for general definitions.
This fact was proved by K.~A.~Muttalib in [Mu]. 
He also argued that biorthogonal ensembles are of certain interest in physics. 

We consider three different cases:
$$
\align
&(1)\quad I=(0,1),\quad \omega(x)=x^{\alpha},\\ 
&(2)\quad I=(0,+\infty),\quad\omega(x)=x^{\alpha}e^{-x},\\  
&(3)\quad I=(-\infty,+\infty),\quad \omega(x)=|x|^{\alpha}e^{-x^2}.
\endalign
$$ 

Naturally, we shall call the second ensemble {\it biorthogonal Laguerre ensemble}, it depends on one real parameter
$\alpha>-1$. 

Ensemble (1) is given by a special case of the Jacobi weight function. Namely, the factor $(1-x)^{\beta}$
of the general Jacobi weight is absent in our formula for $\omega(x)$. Our techniques does not allow to handle the case 
of arbitrary
$\beta$, so we shall set $\beta=0$ but we shall still use the words {\it biorthogonal Jacobi ensemble} in this case.
Thus, our Jacobi biorthogonal ensemble also depends on one real parameter $\alpha>-1$.

The weight function in (3) is, in contrary, more general than the Hermite weight function, which corresponds to $\alpha=0$.
However, we shall deal with this more general case and use the words {\it biorthogonal Hermite ensemble}. 
This ensemble depends on one real parameter $\alpha>-1$.
 
We prove that in all these cases after an appropriate scaling there
 exists a limit of the correlation functions,
 and we explicitly compute the limit kernels thus obtained. 
 
 It turns out that the limit kernel in Jacobi 
 and Laguerre cases is the same and it equals
$$
\Cal K^{(\alpha,\theta)}(x,y)=\sum_{k,l=0}^{\infty}\frac{(-1)^{k}x^{k}}
{k!\Gamma(\frac{\alpha+1+ k}{\theta} )}\frac{(-1)^ly^{\theta l}}
{l!\Gamma({\alpha+1+\theta l})}\frac{\theta}{\alpha+1+ k+\theta l}.
$$
This kernel can be represented as the integral
$$
\Cal K^{(\alpha,\theta)}(x,y)=\theta\int_0^1J_{\frac{\alpha+1}{\theta},\frac{1}\theta}
(xt)\cdot J_{{\alpha+1},{\theta}}((yt)^{\theta})t^{\alpha}dt
$$
where
$$
J_{a,b}(x)=\sum_{m=0}^{\infty}\frac{(-x)^m}{m!\Gamma(a+bm)}
$$
is the Wright's generalized Bessel function, see [Wr], [E2, 18.1].

The limit correlation functions have the form
$$
\rho^{Jac(\alpha,\theta)}_n(x_1,\ldots,x_n)=
\rho^{Lag(\alpha,\theta)}_n(x_1,\ldots,x_n)=
\pro_{i=1}^nx_i^{\alpha}\cdot\det\left[ \Cal K^{(\alpha,\theta)}(x_i,x_j)\right]_{i,j=1}^n.
$$
These functions are defined on $\Bbb R_+$.
  
The limit kernel in the Hermite case is expressed via $\Cal K^{(\alpha,\theta)}
(x,y)$ and has the form
$$
\Cal K^{Her(\alpha,\theta)}(x,y)=\Cal K^{\left(\frac{\alpha-1}{2},\theta\right)}
(x^2,y^2)+x^{\theta}y\cdot 
\Cal K^{\left(\frac{\alpha+\theta}{2},\theta\right)}(x^2,y^2).
$$
 The limit correlation functions are
$$
\rho^{Her(\alpha,\theta)}_n(x_1,\ldots,x_n)=\pro_{i=1}^n|x_i|^{\alpha}\cdot\det
\left[ \Cal K^{Her(\alpha,\theta)}(x_i,x_j)\right]_{i,j=1}^n.
$$
They are defined on the whole real axis. 

Two limit kernels $|xy|^{\alpha/2}\Cal K^{Her(\alpha,\theta)}(x,y)$ and
 $(xy)^{\alpha/2}\Cal K^{(\alpha,\theta)}(x,y)$ can be considered
as biorthogonal generalizations of sine and Bessel kernels, see below 
Examples 5.5 and 3.5, respectively. 

Though these kernels are not symmetric anymore, 
the transposition $x\leftrightarrow y$ leads to a non-trivial 
symmetry in the asymptotics of Laguerre and Hermite ensembles.
 Namely, the change of parameters 
$$
\alpha \mapsto \frac{\alpha+1}\theta -1,\quad \theta\mapsto \frac 1\theta
$$ 
turns out to be equivalent to the transformation $x\mapsto x^{\theta}$ 
of the phase space, see Corollaries 4.7 and 5.6 below.
Finite point ensembles do not possess this symmetry,
it appears only in the asymptotics. 

The method that we use in this paper also provides a new approach to the asymptotics of classic orthogonal polynomial
 ensembles.

The paper is organized as follows. Section 1 is the introduction. In Section 2 we prove some general
statements for further purposes. Section 3 deals with the biorthogonal Jacobi ensemble. In Section 
4 we work with the biorthogonal Laguerre ensemble. The trick that we use there
to obtain a formula for
Christoffel-Darboux kernels originates from our work on  
stochastic point processes arising in the representation theory 
of the infinite symmetric group (see [O], [B], [BO]). We hope to explain 
this connection in subsequent publications. In Section 5 we compute the asymptotics of the biorthogonal Hermite
ensemble. Section 6 is an appendix, it contains rigorous proofs of two main 
theorems. Heuristic proofs of these
theorems can be found at the appropriate places of the main text.  

I am very grateful to G.~I.~Olshanski and A.~A.~Kirillov for numerous and helpful discussions.  
 
\head 2. Generalities \endhead
We start with a discussion of well-known orthogonal polynomial ensembles,
 see [NW1], [Me]. 
Consider an ensemble of $N$ points on a (possibly infinite) interval $(a,b)$ of the real
 line with the joint probability density of the form
$$
p(x_1,\ldots,x_N)=const\cdot\pro_{i=1}^N\omega(x_i)\pro_{i<j}(x_i-x_j)^2
\tag 2.1
$$
where $\omega(x)$ is some positive weight function on $(a,b)$. We can explicitly compute all correlation functions of this ensemble using orthogonal polynomials 
$$\{p_k(x)\}_{k=0}^\infty,\quad
\operatorname{deg}p_k=k;
$$
$$
\inte_a^bp_k(x)p_l(x)\omega(x)dx=\delta_{kl}.
$$ 
Namely, we consider the {\it $N$th Christoffel-Darboux  kernel}
$$
K_N(x,y)=\sum_{i=0}^{N-1}p_i(x)p_i(y).
\tag 2.2
$$

The integral operator in $L
^2([a,b],\omega(x)dx)$ with this kernel is the
orthogonal projection on the $N$--dimensional subspace
$Span\{1,x,\ldots,x^{N-1}\}$.

The {\it correlation functions}
$$
\rho_n(x_1,\ldots,x_n)=\frac{N!}{(N-n)!}\inte_a^b\cdots\inte_a^b p(x_1,\ldots,x_N)dx_{n+1}\cdots dx_N
\tag 2.3
$$
have the following form (see [D], [Me], [NW1])
$$
\rho_n(x_1,\ldots,x_n)=\pro_{i=1}^n \omega(x_i)\cdot \det [K_N(x_i,x_j)]_{i,j=1}^n.
\tag 2.4
$$
The {\it Christoffel-Darboux formula} (see [S], [E1])
$$
K_N(x,y)=\frac{p_N(x)p_{N-1}(y)-p_{N-1}(x)p_{N}(y)}{x-y}
\tag 2.5
$$
provides a convenient way of analyzing the asymptotic behaviour of $K_N(x,y)$ when $N\to \infty$.

Let us look at $K_N(x,y)$ from another point of view. Clearly, the kernel has the form 
$$
K_N(x,y)=\sum_{k,l=0}^{N-1}c_{kl} x^ky^l,
\tag 2.6
$$
where $c_{kl}$ are some constants. There is a simple way of saying what these constants are. Set
$$
g_{ij}=\inte_a^bx^{i+j}\omega(x)dx.
$$
Let 
$$
C=(c_{kl})_{k,l=0}^{N-1},\quad G=(g_{ij})_{i,j=0}^{N-1}.
$$
\proclaim{Proposition 2.1} With the preceding notation
$$
C=G^{-1}.
\tag 2.7
$$
\endproclaim
\demo{Proof}
By definition of the Christoffel-Darboux kernel, for all $j=0,1,\ldots,N-1$
$$
\inte_a^bK_N(x,y)y^jdy=x^j.
$$ 
On the other hand
$$
\inte_a^bK_N(x,y)y^jdy=\inte_a^b\sum_{k,l=0}^{N-1}c_{kl} x^ky^l\cdot y^jdy=\sum_{k,l=0}^{N-1}c_{kl} g_{lj}x^k.
$$
By equating these two expressions, we get
$$
\sum_{l=0}^{N-1}c_{kl} g_{lj}=\delta_{kj}. \qed
$$
\enddemo
Proposition 2.1 and  (2.6) provide another way of analyzing the kernel $K_N(x,y)$. As we shall see, this will also work for biorthogonal ensembles, where Christoffel-Darboux formula becomes rather complicated.  

Let us introduce the {\it biorthogonal ensembles}
in the following very general way. Consider an ensemble of $N$ points on  $(a,b)\subset \Bbb R$ with the joint probability density of the form
$$
p^b(x_1,\ldots,x_N)=const\cdot\pro_{i=1}^N\omega(x_i)\cdot\det[ \xi_i(x_j)]_{i,j=1}^{N}\cdot\det[\eta_i(x_j)]_{i,j=1}^{N},
\tag 2.8
$$
where $\xi_i(x),\eta_i(x);\  i=1,2,\ldots$ are some functions defined on $(a,b)$. The formula (2.1) is clearly a special case of (2.8), take $\xi_i(x)=\eta_i(x)=x^{i-1}$. 

Suppose, we managed to biorthogonalize $\{\xi_i\}$ and $\{\eta_i\}$ with respect to the pairing
$$
\langle\xi,\eta\rangle=\inte_a^b\xi(x)\eta(x)\omega(x)dx.
\tag 2.9
$$
That is, we have two other systems of functions, say $\{\zeta_i(x)\}_{i=1}^{N}$ and $\{\psi_i(x)\}_{i=1}^{N}$, such that
$$
\zeta_i\in Span\{\xi_1,\ldots,\xi_i\},\quad \psi_j\in Span\{\eta_1,\ldots,\eta_j\};
$$
$$
\langle \zeta_i,\psi_j\rangle=\delta_{ij}
$$
for all possible $i$ and $j$.

Exactly the same argument as for orthogonal ensembles proves that the correlation functions (defined as in (2.3)) of a biorthogonal ensemble have the form
$$
\rho^b_n(x_1,\ldots,x_n)=\pro_{i=1}^n \omega(x_i)\cdot \det [K^b_N(x_i,x_j)]_{i,j=1}^n,
\tag 2.10
$$
where
$$
K^b_N(x,y)=\sum_{i=1}^{N}\zeta_i(x)\psi_i(y).
\tag 2.11
$$
is an analog of the Christoffel-Darboux kernel (2.2).

In the next section we shall restrict ourselves to the case
$$
\xi_i(x)=x^{i-1},\quad \eta_i(x)=x^{\theta (i-1)};\ \theta>0.
$$
Then $\{\zeta_i\}$ and $\{\psi_i\}$ are so-called {\it biorthogonal polynomials}, see [K1]. 
In this case determinantal formula (2.10) was proved in [Mu].

(In Section 5, however, it will be more convenient to take
$$
\xi_i(x)=x^{\theta(i-1)},\quad \eta_i(x)=x^{i-1}.
$$
This transposition, clearly, does not change the ensemble (2.8) and its
correlation functions (2.10), but it interchanges the variables $x$
and $y$ in the
Christoffel--Darboux kernel (2.11).)

Unfortunately, there is no simple analog of Christoffel-Darboux formula 
in the biorthogonal case.
For example, for integer $\theta=k\in \Bbb N$, the following formula
is proved in [I] (cf. (2.5))
$$
K^b_N(x,y)=\frac{1}{x^k-y^k}\Biggl
(\sum\limits_{\Sb r,s\geq 0\\ r+s\leq k-1\endSb}
\alpha_{rs}\zeta_{N+s+1}(x)\psi_{N-r}(y)+
\beta\zeta_{N}(x)\psi_{N+1}(y)\Biggr),
$$
where $\alpha_{rs}$, $\beta$ are some constants.  
But even in the simplest case $\theta=2$ it requires a lot of technical work to compute these constants
 for known biorthogonal polynomials. 

We shall use another approach. The kernel $K_N^b(x,y)$ has the following form (cf.  (2.6))
$$
K^b_N(x,y)=\sum_{k,l=1}^{N}c^b_{kl} \xi_k(x)\eta_l(y)
\tag 2.12
$$
for some constants $c^b_{kl}$.
 Set
$$
g^b_{ij}=\inte_a^b\xi_j(x)\eta_i(x)\omega(x)dx.
$$
Let 
$$
C^b=(c^b_{kl})_{k,l=1}^{N},\quad G^b=(g^b_{ij})_{i,j=1}^{N}.
$$
It turns out that we have the exact analog of Proposition 2.1.
\proclaim{Proposition 2.2} With the preceding notation
$$
C^b=(G^b)^{-1}.
\tag 2.13
$$
\endproclaim
    The proof is just the same. 

    To conclude this section, we shall consider even more general situation, where an analog of Propositions 2.1 and 2.2 holds.

Let us fix two (possibly infinite) intervals $(a_1,b_1)$ and
 $(a_2,b_2)$ of the real line. We consider a distribution $\goth p(x_1,\ldots,x_N;y_1,\ldots,y_N)$ defined on $(a_1,b_1)^N\times (a_2,b_2)^N$ of the form
$$
\goth p(x_1,\ldots,x_N;y_1,\ldots,y_N)=const\cdot \det w(x_i,y_j)\cdot  \det \xi_i(x_j)\cdot \det \eta_i(y_j),
\tag 2.14
$$ 
all subscripts vary from $1$ to $N$. Here $w(x,y)$ is a 
(generalized) function on $(a_1,b_1)\times (a_2,b_2)$;
 $\xi_i(x),\eta_i(x)$ are some (generalized)
  functions defined on $(a_1,b_1)$ and $(a_2,b_2)$, respectively. 

We choose the constant so that
$$
\inte_{ (a_1,b_1)^N\times (a_2,b_2)^N}\goth p(x_1,\ldots,x_N;y_1,\ldots,y_N)dx_1\cdots dx_Ndy_1\cdots dy_N=1.
$$

Suppose again that we managed to biorthogonalize the systems $\{\xi_i\}$ and $\{\eta_i\}$ with respect to the pairing
$$
\langle\xi,\eta\rangle=\inte_{a_1}^{b_1}\inte_{a_2}^{b_2}\xi(x)\eta(y)w(x,y)dxdy.
\tag 2.15
$$
Thus, we have somehow constructed two systems $\{\zeta_i(x)\}_{i=1}^{N}$ and $\{\psi_i(y)\}_{i=1}^{N}$, such that
$$
\zeta_i\in Span\{\xi_1,\ldots,\xi_i\},\quad \psi_j\in Span\{\eta_1,\ldots,\eta_j\};
$$
$$
\langle \zeta_i,\psi_j\rangle=\delta_{ij}
$$
for all possible $i$ and $j$.

If $a_1=a_2=a$, $b_1=b_2=b$, and $w(x,y)=\omega(x)\delta(x-y)$, then we
return to the situation described above: (2.15) and (2.9) will coincide.

It turns out that  if we compute the ``correlation functions'' of the
``measure''
$\goth p$ only in $(a_1,b_1)$ (i.e., we integrate $\goth p$ over all $N$ $y$'s and
over some $x$'s), then these correlation functions also have determinantal form. (In fact, as was recently
 proved in [EM], all correlation functions of measures of the type (2.14) have determinantal form, but this statement is much
 harder.)
 
The following statement is proved in [MS]. 
  
\proclaim{Proposition 2.3} With the preceding notation
$$
\align
\frac{N!}{(N-n)!}\inte_{(a_1,b_1)^{N-n}\times (a_2,b_2)^N}
\goth p(x_1,\ldots,x_N;y_1,\ldots,y_N)dx_{n+1}&\cdots dx_Ndy_1\cdots dy_N\\=
&\det [\tilde K_N(x_i,x_j)]_{i,j=1}^n,
\tag 2.16
\endalign
$$
where
$$
\tilde K_N(x,t)=\su_{i=1}^N\zeta_i(x)\inte_{a_2}^{b_2}\psi_i(y) w(t,y)dy.
\tag 2.17
$$
\endproclaim
We shall use Proposition 2.3 in Section 4.

Note that if $a_1=a_2=a$, $b_1=b_2=b$, and $w(x,y)=\omega(x)\delta(x-y)$, then
$$
\tilde K_N(x,t)=\omega(t)K^b_N(x,t).
$$

As above, we can write 
$$
\su_{i=1}^N\zeta_i(x)\psi_i(y)=\sum_{k,l=1}^{N}\tilde c_{kl} \xi_k(x)\eta_l(y)
\tag 2.18
$$
for some constants $\tilde c_{kl}$. 

Set
$$
\tilde g_{ij}=\inte_{a_1}^{b_1}\inte_{a_2}^{b_2}\xi_j(x)\eta_i(y)w(x,y)dxdy;
$$
$$
\tilde C=(\tilde c_{kl})_{k,l=1}^{N},\quad \tilde G=(\tilde g_{ij})_{i,j=1}^{N}.
$$ 
Quite similarly to Propositions 2.1 and 2.2 we get the following assertion.
\proclaim{Proposition 2.4} With the preceding notation
$$
\tilde C=\tilde G^{-1}.
\tag 2.19
$$
\endproclaim
By (2.17) and (2.18) we get 
$$
\tilde K_N(x,t)=\sum_{k,l=1}^{N}\tilde c_{kl} \xi_k(x)\inte_{a_2}^{b_2}\eta_l(y)w(t,y)dy.
\tag 2.20
$$

Note that everywhere above we {\it supposed} that the
corresponding biorthogonal systems exist. The following simple statement
 will guarantee the existence in all our further examples.
\proclaim{Proposition 2.5} If all principal minors of $\tilde G$ are not
 zero, then there exist\linebreak biorthogonal systems 
 $\{\zeta_i(x)\}_{i=1}^{N}$ 
 and $\{\psi_i(y)\}_{i=1}^{N}$.
\endproclaim
\demo{Proof} The hypothesis implies that $\tilde G$ possesses
 a Gauss decomposition:
 it can be represented as the product of a lower triangular and an upper
  triangular matrices. 
 Thus, there exist a lower triangular matrix $L=(l_{ij})$ and an upper
 triangular matrix $U=(u_{ij})$ such that
$$
L{\tilde G} U=\operatorname{Id}.
$$
Set
$$
\zeta_i=\su_{j=1}^Nu_{ji}\xi_j;\qquad \psi_j=\su_{i=1}^Nl_{ji}\eta_i.
$$
A straightforward check shows that these systems are biorthogonal. \qed
\enddemo
\head 3. Biorthogonal Jacobi ensemble \endhead
Our goal in this section is to study the asymptotic 
behaviour of the $N$-point ensemble on $(0,1)$ with the joint
 probability density (cf. (2.8)) 
$$
\aligned
p^{Jac}_N(x_1,\ldots,x_N)=const\cdot \pro_{i=1}^Nx_i^\alpha\pro_{i<j}&\left[(x_i-x_j)(x_i^{\theta}-x_j^{\theta})\right]\\=const\cdot \pro_{i=1}^Nx_i^\alpha&\cdot \det x_i^{j-1}\cdot\det x_i^{\theta (j-1)}
\endaligned
\tag 3.1
$$
where
$\theta>0$ and $\alpha>-1$. We call this ensemble the {\it biorthogonal 
Jacobi ensemble}.  Let us compute the matrix of pairwise scalar products
 $(g_{ij}^{Jac})$. We have
$$
g_{ij}^{Jac}=\inte_0^1x^{j-1+\theta (i-1)}x^{\alpha}dx=\frac{1}{j+\theta 
(i-1)+\alpha}.
\tag 3.2
$$
To invert this matrix we shall use the following lemma.
\proclaim{Lemma 3.1} Let $A=\{A_1,\ldots,A_N\}$ and
 $B=\{B_1,\ldots,B_N\}$ be two sequences of complex numbers
  such that $A_i+B_j\neq 0$ for all $i,j=1,\ldots,N$, and
  $A_i\neq A_j,\ B_i\neq B_j$ for $i\neq j$. Then
$$
\left(\frac{1}{A_i+B_j}\right)^{-1}=(C_{kl})
$$
where
$$
C_{kl}=\frac
{\prod_{i=1}^N\left[(B_i+A_l)(A_i+B_k)\right]}
{\prod_{i\neq l}(A_l-A_i)
\prod_{j\neq k}(B_k-B_j)}
\frac{1}{A_l+B_k}.
\tag 3.3
$$
\endproclaim
\demo{Proof} As is  known, the elements of the inverse matrix are the
 ratios of the cofactors of
  corresponding elements of the initial matrix and the determinant 
  of the initial matrix. The determinant of our matrix 
$$
M=\left(\frac{1}{A_i+B_j}\right)
$$
 is well-known {\it Cauchy determinant}, see [W]:
$$
\det M=\frac{\pro_{i<j}\left[(A_i-A_j)(B_i-B_j)\right]}{\pro_{i=1}^N\pro_{j=1}^N(A_i+B_j)}.
$$
Every submatrix of $M$  has the same form as $M$ itself for some other
 sets $A$ and $B$. Then we can use the formula for the Cauchy determinant
  for computing any minor of   $M$. In particular, we can compute all cofactors,
   and, thus, the inverse matrix. The result is exactly (3.3). \qed
\enddemo
\proclaim{Proposition 3.2} The inverse of the Gram matrix $(g_{ij}^{Jac})$
 has the form
$$
(g_{ij}^{Jac})^{-1}=\theta\left(\frac{\left(\frac{k+\alpha}{\theta}\right)_N}
{(k-1)!(N-k)!}\cdot\frac{(\theta (l-1)+\alpha+1)_N}{(l-1)!(N-l)!}\cdot
 \frac{(-1)^{k+l}}{k+\theta (l-1)+\alpha}\right)_{k,l=1}^N
\tag 3.4
$$
where $(a)_m=a(a+1)\cdots(a+m-1)$ stands for the Pohgammer symbol.
\endproclaim
\demo{Proof} Direct application of Lemma 3.1 for 
$$
A_i=\theta(i-1),\quad B_i=i+\alpha;\quad i=1,\ldots, N.\qed
$$
\enddemo
\proclaim{Proposition 3.3} The correlation functions of the biorthogonal 
$N$-point Jacobi ensemble have the form
$$
\rho_{nN}^{Jac}(x_1,\ldots,x_n)=\pro_{i=1}^nx_i^{\alpha}\cdot \det[K_N^{Jac}(x_i,x_j)]_{i,j=1}^n
$$
where
$$
K_N^{Jac}(x,y)=\theta\su_{k,l=1}^N
\frac{\left(\frac{k+\alpha}{\theta}\right)_Nx^{k-1}}{(k-1)!(N-k)!}
\cdot\frac{(\theta (l-1)+\alpha+1)_Ny^{\theta(l-1)}}{(l-1)!(N-l)!}
\cdot \frac{(-1)^{k+l}}{k+\theta (l-1)+\alpha}
\tag 3.5
$$
\endproclaim
\demo{Proof} The claim follows from Proposition 2.2 and Proposition 3.2.
 The existence of the Christoffel-Darboux kernel $K_N^{Jac}(x,y)$ is
  guaranteed by Proposition 2.5 because all minors of the matrix
  $(g_{ij}^{Jac})$ are nonzero. \qed
\enddemo
Now we are in a position to compute the asymptotics of our ensemble as $N\to \infty$.

 We shall employ the following 
entire function introduced by E.~M.~Wright, see [Wr], [E2, 18.1(27)]
$$
J_{a,b}(x)=\sum_{m=0}^{\infty}\frac{(-x)^m}{m!\Gamma(a+bm)}
$$
(our notation differs from that used in [Wr], [E2]).
It is closely related to Mittag--Leffler type functions, see [E2, 18.1].

Note that
$$
x^{\frac{a}2}J_{a+1,1}=\sum_{m=0}^{\infty}\frac{(-1)^m x^{m+\frac a2}}
{m!\Gamma(a+m+1)}=J_a(2\sqrt{x})
$$
where $J_a(x)$ is the Bessel function. 

\proclaim{Theorem 3.4} For all $n=1,2,\ldots$ there exists the limit
$$
\lim_{N\to \infty}\frac{1}{N^{n(1+\frac1\theta)}}
\rho^{Jac}_{nN}\left(\frac{x_1}{N^{1+\frac1\theta}},\ldots,
\frac{x_n}{N^{1+\frac1\theta}}\right)=\pro_{i=1}^nx_i^{\alpha}
\cdot\det\left[ \Cal K^{(\alpha,\theta)}(x_i,x_j)\right]_{i,j=1}^n
$$ 
where
$$
\aligned
\Cal K^{(\alpha,\theta)}(x,y)=\sum_{k,l=0}^{\infty}&\frac{(-1)^{k}x^{k}}
{k!\Gamma(\frac{\alpha+1+ k}{\theta} )}\frac{(-1)^ly^{\theta l}}{l!\Gamma({\alpha+1+\theta l})}\frac{\theta}{\alpha+1+ k+\theta l}\\
&=\theta\int_0^1J_{\frac{\alpha+1}{\theta},\frac{1}\theta}(xt)\cdot J_{{\alpha+1},{\theta}}((yt)^{\theta})t^{\alpha}dt.
\endaligned
\tag 3.6
$$
\endproclaim
\demo{Heuristic proof}
If we substitute the asymptotic relations
$$
\frac{\left(\frac{k+\alpha}{\theta}\right)_N}{(N-k)!}\sim \frac{N^{\frac{(k-1)(\theta+1)+\alpha+1}\theta}}{\Gamma\left(\frac{k+\alpha}{\theta}\right)}$$$$
\frac{(\theta (l-1)+\alpha)_N}{(N-l)!}\sim\frac{N^{(l-1)(\theta+1)+\alpha+1}}{\Gamma(\alpha+1+\theta (l-1))}
$$
in (3.5) and shift the summation indices $k$ and $l$ by one, then we get
(3.6).\qed
\enddemo
A rigorous proof of Theorem 3.4 (which is nothing more than a formalization of
the argument above) can be found in Section 6.
\example{Example 3.5} For $\theta=1$, Theorem 3.4 (and Theorem 4.5, see below)
is well known, see [NW2], [F].
 In this case it is easy to check that
$$
(xy)^{\alpha/2}\Cal K^{(\alpha,\theta)}(x,y)=\frac{\varphi_1(x)\varphi_2(y)-\varphi_1(y)\varphi_2(x)}{x-y}
$$
where
$$
\varphi_1(x)=J_{\alpha}(2\sqrt{x}),\quad \varphi_2(x)=x\varphi_1'(x),
$$
is the {\it Bessel kernel}.
\endexample
Let us draw one more interesting corollary from Lemma 3.1.

 Note that in the notation of Section 2 if we know the Christoffel-Darboux kernels (2.11) then we can easily obtain explicit formulas for biorthogonal functions $\{\zeta_i(x)\}$  and $\{\psi_i(x)\}$. Namely, we have the relation
$$
\zeta_m(x)\psi_m(y)=K_m^b(x,y)-K_{m-1}^b(x,y).
$$
By Lemma 3.1 we can compute the kernel $K_m^b(x,y)$ for the systems
$$
\xi_i(x)=x^{a_i},\quad \eta_i(x)=x^{b_i};
$$
$$ a_i+b_j>-1,\quad i,j=1,2,\ldots
$$
in $L^2([0,1],dx)$.
The result of computation of the corresponding biorthogonal systems is expressed in the next statement.
\proclaim{Proposition 3.6} Let $a_1,a_2,\ldots$ and $b_1,b_2,\ldots$ be
two sequences of complex numbers such that
$$ a_i+b_j>-1,\quad i,j=1,2,\ldots
$$
and $a_i\neq a_j,\ b_i\neq b_j$ for $i\neq j$.
Then two systems of functions 
$$
\zeta_n(x)=\sqrt{a_{n}+b_{n}+1}\sum\limits_{i=1}^{n}\frac{\prod_{k=1}^{n-1}(a_i+b_k+1)}{\prod_{k=1,\neq i}^{n}(a_i-a_k)}x^{a_i},\quad n=1,2,\ldots
$$
and
$$
\psi_n(x)=\sqrt{a_{n}+b_{n}+1}\sum\limits_{i=1}^{n}\frac{\prod_{k=1}^{n-1}(b_i+a_k+1)}{\prod_{k=1,\neq i}^{n}(b_i-b_k)}x^{b_i},\quad n=1,2,\ldots
$$
are biorthonormal in $L^2([0,1],dx)$. In other words
$$
\inte_0^1\zeta_m(x)\psi_n(x)dx=\delta_{mn}.
$$
\endproclaim
\example{Example 3.7} For $$a_i=b_i=i-1$$ our biorthogonal systems are classic Jacobi polynomials. 
If $$a_i=i-1,\quad b_i=\theta (i-1)$$ for some $\theta>0$ 
we get explicit formulas for biorthogonal Jacobi polynomials, cf. [MT1].
\endexample 
\head 4. Biorthogonal Laguerre ensemble \endhead
In this section we study the asymptotic behaviour of the $N$-point ensemble
 on $(0,+\infty)$ with the joint probability density
$$
\aligned
p^{Lag}_N(x_1,\ldots,x_N)=const\cdot \pro_{i<j}&\left[(x_i-x_j)(x_i^{\theta}-x_j^{\theta})\right]\cdot\pro_{i=1}^Nx_i^{\alpha}\cdot e^{-x_1-\ldots-x_N}\\=const&\cdot \det x_j^{i-1}\cdot\det x_j^{\alpha+\theta (i-1)}\cdot e^{-x_1-\ldots-x_N}
\endaligned
\tag 4.1
$$
where $\alpha>-1$ and $\theta>0$.

We call it the $N$-point {\it biorthogonal Laguerre ensemble}. It turns out 
that the asymptotics of this ensemble is governed by the same kernel as that 
of biorthogonal Jacobi ensemble, see the previous section. However, in this
case it is much more difficult to show. Here we do not know an analog of
Lemma 3.1 for inverting the matrix of scalar products. However, the 
following trick allows to reduce the computation of the
correlation functions to Lemma 3.1. 

Let us introduce $N$ new variables $y_1,\ldots,y_N$ and a new distribution in 
$2N$ variables
$$
\goth p_N(x_1,\ldots,x_N;y_1,\ldots,y_N)=p^{Lag}_N(x_1,\ldots,x_N)
\pro_{i=1}^N\delta(y_i).
$$
It is quite clear that
$$
\multline
\rho_{nN}^{Lag}(x_1,\ldots,x_n)=\frac{N!}{(N-n)!}
\int p_N^{Lag}(x_,\ldots,x_N)dx_{n+1}\cdots dx_N\\=
\frac{N!}{(N-n)!}
\int\goth p_N(x_1,\ldots,x_N;y_1,\ldots,y_N)dx_{n+1}\cdots dx_Ndy_{1}\cdots dy_N.
\endmultline
$$
Thus, it suffices to compute the correlation functions of $\goth p_N$. To do
this we shall use Proposition 2.3, but, first, we need to show that 
$\goth p_N$ has the form (2.14).

\proclaim{Proposition 4.1} With the preceding notation
$$
\goth p_N(x_1,\ldots,x_N;y_1,\ldots,y_N)=const\cdot \det w(x_i,y_j)\cdot 
 \det \xi_i(x_j)\cdot \det \eta_i(y_j),
$$
all indices vary from $1$ to $N$, and
$$
\xi_i(x)=x^{\alpha+ N+\theta(i-1)},\quad \eta_i(y)=\delta^{(i-1)}(y)
$$
for all $i=1,\ldots,N$;
$$
w(x,y)=\frac{e^{-x-y}}{x+y}.
$$
\endproclaim
\demo{Proof} Using the formula for Cauchy determinant mentioned in the proof
of Lemma 3.1, we get
$$
\det w(x_i,y_j)=\frac{\pro_{i<j}\left[(x_i-x_j)(y_i-y_j)\right]}{\pro_{i,j}(x_i+y_j)}\cdot e^{-\su_{i=1}^N(x_i+y_i)}
$$
Next,
$$
\pro_{i<j}(y_i-y_j)\cdot \det \delta^{(i-1)}(y_j)=(-1)^{\frac{n(n-1)}2}N!\pro_{i=1}^N\delta(y_i)
$$
because $y^m\cdot\delta^{(n)}(y)=\delta_{mn}$ for $m\geq n$.
Finally,
$$
\multline
\goth p_N(x_1,\ldots,x_N;y_1,\ldots,y_N)\\
=const\cdot(-1)^{\frac{n(n-1)}2}\det x_j^{\alpha+N+\theta(i-1)}\frac{\pro_{i<j}(x_i-x_j)}{\pro_{i,j}(x_i+y_j)}\cdot e^{-\su_{i=1}^N(x_i+y_i)}\cdot N!\pro_{i=1}^N\delta(y_i)\\
=const\cdot \pro_{i<j}\left[(x_i-x_j)(x_i^{\theta}-x_j^{\theta})\right]\cdot \pro_{i=1}^Nx_i^{\alpha}\cdot e^{-x_1-\ldots-x_N}\cdot\pro_{i=1}^N\delta(y_i).\qed
\endmultline
$$
\enddemo
Now, according to Section 2, we need to compute
$$
\tilde g_{ij}=\inte_{0}^{+\infty}\inte_{0}^{+\infty}\xi_j(x)\eta_i(y)w(x,y)dxdy.
$$
\proclaim{Proposition 4.2}
$$
\tilde g_{ij}=\inte_0^{+\infty}\inte_0^{+\infty}\frac{x^{\alpha+N+\theta(j-1)}
\delta^{(i-1)}(y)e^{-x-y}}{x+y}dxdy=
\frac{\Gamma(1+\alpha+N+\theta(j-1))}{\alpha+N+\theta(j-1)-(i-1)}
$$
\endproclaim
Before proving Proposition 4.2, let us explain our achievement. By Lemma 3.1 we can now easily compute the inverse matrix $$\tilde C=(\tilde c_{kl})=(\tilde g_{ij})^{-1}.$$
\proclaim{Corollary 4.3} With the preceding notation
$$
\tilde c_{kl}=\theta\,\frac{\Gamma(1+\alpha+\theta(k-1))^{-1}}{(k-1)!(N-k)!}\,
\frac{\left(\frac{\alpha+N-(l-1)}{\theta}\right)_N}{(l-1)!(N-l)!}\,
\frac{(-1)^{N+k+l+1}}{\alpha+N+\theta(k-1)-(l-1)}
\tag 4.2
$$
\endproclaim
\demo{Proof of Corollary 4.3} Apply Lemma 3.1 for
$$
A_i=-(i-1),\quad B_i=\alpha+N+\theta(i-1),\quad i=1,\ldots,N. \qed
$$
\enddemo
\demo{Proof of Proposition 4.2}
It is easy to check that
$$
\inte_0^{+\infty}\inte_0^{+\infty}\frac{x^a}{\Gamma(a+1)}\frac{y^b}{\Gamma(b+1)}\frac{e^{-x-y}}{x+y}dxdy=\frac{1}{a+b+1}
\tag 4.3
$$
if, say, $\Re a$ and $\Re b >0$. (Change of variables $r=x+y$, $s=\frac x{x+y}$
reduces the integral to the product of Euler's gamma and beta integrals.)

 As is well known, there exists a distribution 
$$
\phi_c(u)=\frac{u_+^{c}}{\Gamma(c+1)}
$$
which depends on $c$ analytically, 
such that for $c>-1$ it is just an integrable function ${u^{c}}/{\Gamma(c+1)}$ for $u>0$ and $0$ for $u\leq 0$. For $c\leq -1$ it is defined via analytic continuation. In particular, we always have the relation $\phi_c'=\phi_{c-1}$ and, thus, for any positive integer $k$
$$
\phi_{-k}(u)=\delta^{(k-1)}(u).
$$
Analytic continuation of (4.3) gives
the relation
$$
\inte_0^{+\infty}\inte_0^{+\infty}\frac{x^a}{\Gamma(a+1)}\phi_b(y)\frac{e^{-x-y}}{x+y}dxdy=\frac{1}{a+b+1}
$$
for, say, $\Re a>0$ and $\Re (a+b+1)>0$. Our claim is the special case of this formula for
$$
a=\alpha+N+\theta(i-1),\quad b=-j.\qed
$$
\enddemo

Thus, we have inverted the matrix $(\tilde g_{ij})$, and as the result we get
a formula for the correlation functions.
\proclaim{Theorem 4.4} The correlation functions of the $N$-point biorthogonal Laguerre ensemble have the form
$$
\rho_{nN}^{Lag}(x_1,\ldots,x_n)=\pro_{i=1}^nx_i^{\alpha}e^{-x_i}\cdot 
\det[K_N^{Lag}(x_i,x_j)]_{i,j=1}^n
$$
where
$$
\aligned
K_N^{Lag}(x,y)=
\theta \sum_{k,i=0}^{N-1}\sum_{r=i}^{N-1}
&\frac{\Gamma\left(N+\frac{i+\alpha+1}{\theta}\right)}
{\Gamma(\alpha+\theta k +1)\Gamma\left(\frac{i+\alpha+1}{\theta}\right)}
\\ &\times \frac{(-1)^{i+k}}{k!(N-k-1)!i!(r-i)!}\frac{x^{\theta k}y^{r}}
{\alpha+\theta k +i+1}.
\endaligned
\tag 4.4
$$
\endproclaim
\demo{Proof} By the definition of $\goth p_N$,
$$
\multline
\rho_{nN}^{Lag}(x_1,\ldots,x_n)=\frac{N!}{(N-n)!}\int p_N^{Lag}(x_,\ldots,x_N)dx_{n+1}\cdots dx_N\\=
\frac{N!}{(N-n)!}\int\goth p_N(x_1,\ldots,x_N;y_1,\ldots,y_N)dx_{n+1}\cdots dx_Ndy_{1}\cdots dy_N.
\endmultline
$$
Thus, we can apply Propositions 2.3 and 2.4. We have
$$
\multline
\int \psi_l(y)w(t,y)dy=\int\delta^{(l-1)}(y)\frac{e^{-t-y}}{t+y}dy\\=(-1)^{l-1}\frac{\partial^{l-1}}{\partial y^{l-1}}\frac{e^{-t-y}}{t+y}\biggm|_{y=0}=\su_{s=0}^{l-1}\frac{(l-1)!}{(l-1-s)!}t^{-s-1}e^{-t}.
\endmultline
$$
Then, by Proposition 2.3, we get the determinantal formula
$$
\rho_{nN}^{Lag}(x_1,\ldots,x_n)=\det \tilde K_N(x_i,x_j)
$$
where, see (4.2),
$$
\multline
\tilde K_N(x,t)=\su_{k,l=1}^N\tilde c_{kl}\su_{s=0}^{l-1}\frac{(l-1)!}{(l-1-s)!}t^{-s-1}e^{-t}=
\theta x^{\alpha}e^{-t}\frac{x^N}{t^N}\\ \times \su_{k,l=1}^N\su_{s=0}^{l-1}\frac{\Gamma(1+\alpha+\theta(k-1))^{-1}}{(k-1)!(N-k)!}\frac{\left(\frac{\alpha+N-(l-1)}{\theta}\right)_N}{(l-1-s)!(N-l)!}\\ \times\frac{(-1)^{N+k+l+1}x^{\theta(k-1)}t^{N-s-1}}{\alpha+N+\theta(k-1)-(l-1)}.
\endmultline
$$
Note that the factor $(x/t)^N$ disappears, when we take the determinant,
 and the factor $x^{\alpha}e^{-t}$ produces the factor $\pro_{i=1}^nx_i^{\alpha}e^{-x_i}$ outside the determinant.
Introducing new summation indices 
$$
i=N-l,\quad r=N-s-1
$$
and shifting index $k$ by one, we arrive at our claim.\qed
\enddemo

One may ask whether the kernel $K_N^{Lag}(x,y)$ defined above is of the
Christoffel--Darboux type in the sense of Section 2 (see (2.11)). The answer
is positive: if we take (using the notation of Section 2)
$$
\xi_i(x)=x^{\theta(i-1)},\quad \eta_i(x)=x^{i-1}
$$
then
$$
K_N^{Lag}(x,y)=\su_{i=1}^N \zeta_i(x)\psi_i(y),
$$
see Proposition 5.1.

Now we are ready to pass to the limit $N\to \infty$.
\proclaim{Theorem 4.5} For all $n=1,2,\ldots$ there exists the limit
$$
\lim_{N\to \infty}\frac{1}{N^{\frac n\theta}}
\rho^{Lag}_{nN}\left(\frac{x_1}{N^{\frac 1\theta}},
\ldots,\frac{x_n}{N^{\frac1\theta}}\right)=\pro_{i=1}^nx_i^{\alpha}
\cdot\det\left[ \Cal K^{(\alpha,\theta)}(x_i,x_j)\right]_{i,j=1}^n
$$ 
where
$\Cal K^{(\alpha,\theta)}(x,y)$
is defined in (3.6).
\endproclaim
\demo{Heuristic proof}Indeed, if we substitute in (4.4) the asymptotic relation
$$
\frac{\Gamma\left(n+\frac{i+\alpha+1}{\theta}\right)}
{(N-k-1)!}\sim N^{\frac{i+\alpha+1}{\theta}+k}
$$
then we get
$$
\multline
\frac 1{N^{\frac1\theta}}\left(\frac x{N^{\frac1\theta}}\right)^\alpha
 e^{-{x}/{N^{\frac1\theta}}} K_N^{Lag}\left(\frac{x}{N^{\frac1\theta}},
 \frac y{N^{\frac1\theta}}\right)\sim x^{\alpha}e^{-{x}/{N^{\frac1\theta}}}
  \theta\\ \times\sum_{k,i=0}^{n-1}
\frac{(-1)^{i+k}x^{\theta k}y^i}{k!\Gamma(\alpha+\theta k +1)i!
\Gamma\left(\frac{i+\alpha+1}{\theta}\right)(\alpha+\theta k +i+1)}
\sum_{r=i}^{n-1}\frac{(y/N^{\frac1\theta})^{r-i}}{(r-i)!}.
\endmultline
$$
But
$$
\sum_{r=i}^{n-1}\frac{(y/N^{\frac1\theta})^{r-i}}{(r-i)!}=1+o(1),\quad 
e^{-{x}/{N^{\frac1\theta}}}=1+o(1),
$$
 and we arrive at our asymptotic kernel.\qed 
\enddemo
A justification of the heuristic argument above is given
in the appendix (Section 6).
\example{Remark 4.6} Asymptotic behaviour of biorthogonal Jacobi and Laguerre 
ensembles is the same, while the scaling factors are different.
\endexample

Now let us point out a certain symmetry in the kernel $\Cal
K^{(\alpha,\theta)}(x,y)$ as defined by (3.6):
$$
\Cal K^{(\alpha,\theta)}(y^{\frac 1\theta},x^{\frac 1\theta})=\frac 1\theta\,
\Cal K^{\left(\frac{\alpha+1}\theta -1,
\frac 1\theta\right)}
(x,y).
$$
It follows that the correlation functions (defined by the
determinantal formula of Theorem 4.5), or rather the correlation {\it
measures}
$$
\prod\limits_{i=1}^nx_i^{\alpha}
\cdot\det\left[ \Cal K^{(\alpha,\theta)}(x_i,x_j)\right]_{i,j=1}^n
\cdot dx_1\dots dx_n
$$
are stable under the change of parameters
$$
\alpha \mapsto \frac{\alpha+1}\theta -1,\quad \theta\mapsto \frac 1\theta
\tag 4.5
$$ 
combined with the transformation of the phase space 
$x\mapsto x^{\frac1\theta}, \ y\mapsto y^{\frac1\theta}$.

For the Jacobi ensemble this is perfectly understandable, the transformation
 of the phase space is equivalent to the change of parameters, because if
  $x_i\mapsto x_i^{\frac1\theta}$ then
$$
\multline
const\cdot \pro_{i=1}^Nx_i^\alpha\pro_{i<j}\left[(x_i-x_j)
(x_i^{\theta}-x_j^{\theta})\right]dx_1\cdots dx_N\\
\longmapsto const\cdot 
\pro_{i=1}^Nx_i^{\frac{\alpha+1}{\theta}-1}
\pro_{i<j}\left[(x_i-x_j)(x_i^{\frac 1\theta}-x_j^{\frac1\theta})\right]dx_1\cdots dx_N
\endmultline
$$
 However, for the Laguerre ensemble this symmetry is rather surprising,
  the change of parameters is not equivalent to a transformation of the phase
   space, because the factor $e^{-x_1-\ldots-x_N}$ does not behave properly. 
   Thus, we get a non-trivial conclusion, which deserves a separate statement.
\proclaim{Corollary 4.7} The asymptotic behaviour of the Laguerre biorthogonal
 ensemble at the left edge of spectrum is invariant with
  respect to the change of parameters (4.5) and the 
 transformation $x\mapsto x^{\frac1\theta}$ of the phase space $(0,+\infty)$.
\endproclaim
\head 5. Biorthogonal Hermite ensemble \endhead
Everywhere below the symbol $x^{\theta}$ stands for $\operatorname{sign}(x)|x|^{\theta}$, and for all integers $k$
$$
x^{\theta k}=(x^\theta)^k.
$$

In this section we are dealing with ensembles on $(-\infty,+\infty)$ with the joint probability densities
$$
\aligned
p^{Her}_N(x_1,\ldots,x_N)=const\cdot \pro_{i=1}^N|x_i|^\alpha \pro_{i<j}\left[(x_i-x_j)(x_i^{\theta}-x_j^{\theta})\right]e^{-x_1^2-\ldots-x_N^2}&\\=const\cdot \pro_{i=1}^N|x_i|^\alpha e^{-x_i^2}\cdot \det x_i^{j-1}\cdot\det x_i^{\theta (j-1)}&
\endaligned
\tag 5.1
$$
where
$\theta>0$ and $\alpha>-1$. We call these ensembles the {\it biorthogonal Hermite ensembles}.

We shall reduce the study of their asymptotics to that of the Laguerre ensembles, see the previous section. To do this we need to introduce a notation for biorthogonal Laguerre and Hermite polynomials. 

Let us denote two sequences of biorthogonal Laguerre polynomials by $\{Z_n^{\alpha}(x,\theta)\}$ and  $\{Y_n^{\alpha}(y,\theta)\}$. That is 
$$
\deg Z_n^{\alpha}(x,\theta)=\deg Y_n^{\alpha}(x,\theta)=n
$$
and 
$$
\int_{0}^{\infty}Z_m^{\alpha}(x^{\theta},\theta)Y_n^{\alpha}(x,\theta)x^{\alpha}e^{-x}dx=\delta_{mn}.
$$

   Such polynomials were explicitly constructed in [K2], [C]:
 $$
 Z_n^{\alpha}(x,\theta)=\su_{j=0}^{n} \binom nj \frac{(-1)^jx^{ j}}
 {\Gamma(\theta j+\alpha+1)},
 $$
 $$
 Y_n^{\alpha}(x,\theta)=\frac 1{n!}\su_{r=0}^n\frac {x^r}{r!}
 \su_{i=0}^r (-1)^i\binom ri \left(\frac {i+\alpha+1}\theta\right)_n.
 $$
  Note that our notation is slightly 
   different from the conventional one, usually $Z_n^{\alpha}(x,\theta)$ is multiplied by 
   $\frac{\Gamma(\theta n+\alpha+1)}{n!}$, but we want our polynomials to be orthonormal. 

Following [MT2] we can construct biorthogonal Hermite polynomials  $\{S^{\alpha}_{n}(x,\theta)\}$  and $\{T^{\alpha}_{n}(y,\theta)\}$ satisfying
$$
\deg S_n^{\alpha}(x,\theta)=\deg T_n^{\alpha}(y,\theta)=n
$$
$$
\int_{-\infty}^{+\infty}S_m^{\alpha}(x^{\theta},\theta)T_n^{\alpha}(x,\theta)|x|^{\alpha}e^{-x^2}dx=\delta_{mn}.
$$
Namely,
$$
\align
&S^{\alpha}_{2n}(x,\theta)=Z_n^{{(\alpha-1)}/{2}}(x^{2},\theta);\\
&T^{\alpha}_{2n}(x,\theta)=Y_n^{{(\alpha-1)}/{2}}(x^{2},\theta);\\
&S^{\alpha}_{2n+1}(x,\theta)=xZ_n^{{(\alpha+\theta)}/{2}}(x^{2},\theta);\\
&T^{\alpha}_{2n+1}(x,\theta)=xY_n^{({\alpha+\theta})/{2}}(x^{2},\theta).
\endalign
$$
Again, our notation differs from the usual one by scalar factors.

It is quite clear that the Christoffel-Darboux kernel for Hermite case can be expressed 
via that for Laguerre case, we shall do this in Proposition 5.3.

\proclaim{Proposition 5.1} With the preceding notation
$$
K_N^{Lag}(x,y)=
\su_{i=0}^{N-1}Z_i^{\alpha}(x^{\theta},\theta)Y_i^{\alpha}(y,\theta).
$$
where $K_N^{Lag}(x,y)$ is defined by (4.4).
\endproclaim
\demo{Idea of the proof}
One can prove this statement by direct verifying that 
$$
 K_{N+1}^{Lag}(x,y)-K_N^{Lag}(x,y)
$$
is equal to
$$
Z_N^{\alpha}(x^{\theta},\theta)Y_N^{\alpha}(y,\theta),
$$
we have explicit formulas for both these expressions. However,
this check is rather
tedious.
\qed 
\enddemo
\example{Remark 5.2}
Proposition 5.1 proves that $K_N^{Lag}(x,y)$ is the Christoffel--Darboux
kernel (2.11) for
$$
\xi_i(x)=x^{\theta(i-1)},\quad \eta_i(x)=x^{i-1}.
$$
Note that this fact and (2.10) provide another
proof of Theorem 4.4. However, it is just a {\it checking} proof, the
explicit
formulas for biorthogonal Laguerre polynomials above do not prompt
a suitable expression for the Christoffel--Darboux kernel (4.4).

In Section 4, using the trick with $\goth p_N$, we, actually, managed to 
derive 
the formula (4.4), and its relative 
simplicity allowed us to analyze the asymptotic
behaviour of the biorthogonal Laguerre ensemble (Theorem 4.5).

Moreover, we have a conceptual proof of the following statement
which generalizes
Proposition 5.1: let $\xi_i(x)=x^{a_i}$ for arbitrary
complex numbers $a_i$, $\Re a_i>0$, and $\eta_i(x)=x^{i-1}$;
 then the trick with 
$\goth p_N$ described in
Section 4 always produces the Christoffel--Darboux type kernel. But the
proof is based on a certain formalism which exceeds the limits of the
present paper.
 
\endexample

Set
$$
K_N^{Her}(x,y)=\su_{i=0}^{N-1}S_i^{\alpha}(x^{\theta},\theta)T_i^{\alpha}(y,\theta).
$$
By (2.10),
$$
\rho_{nN}^{Her}(x_1,\ldots,x_n)=\pro_{i=1}^n|x_i|^\alpha e^{-x_i^2}\cdot \det \left[K_N^{Her}(x_i,x_j)\right]_{i,j=1}^n.
$$
 Let us use a more detailed notation for Christoffel-Darboux kernels and write
 $$ K_N^{Lag(\alpha)}(x,y),\quad K_N^{Her(\alpha)}(x,y)$$ instead of $$K_N^{Lag}(x,y),\quad K_N^{Her}(x,y)$$ indicating the dependence on $\alpha$.
\proclaim{Proposition 5.3} The Christoffel-Darboux kernel for the $N$-point biorthogonal Hermite ensemble has the form
$$
K_N^{Her(\alpha)}(x,y)=K_M^{Lag\left(\frac{\alpha-1}2\right)}(x^2,y^2)+x^{\theta}yK_M^{Lag\left(\frac{\alpha+\theta}2\right)}(x^2,y^2)\quad\text{for}\ N=2M
$$
$$
K_N^{Her(\alpha)}(x,y)=K_M^{Lag\left(\frac{\alpha-1}2\right)}(x^2,y^2)+x^{\theta}yK_M^{Lag\left(\frac{\alpha+\theta}2\right)}(x^2,y^2)\quad\text{for}\ N=2M+1
$$
where $K_N^{Lag(\alpha)}(x,y)$ is defined by (4.4).
\endproclaim
\demo{Proof} Immediately follows from explicit expressions for biorthogonal Hermite polynomials via biorthogonal Laguerre polynomials.\qed
\enddemo
Now we can express the asymptotics of the Hermite ensemble via that of the Laguerre ensemble.
\proclaim{Theorem 5.4} For all $n=1,2,\ldots$ there exists the limit (set $M=\frac N2$)
$$
\lim_{N\to \infty}\frac{1}{M^{\frac n{2\theta}}}\rho^{Her}_{nN}\left(\frac{x_1}{M^{\frac1{2\theta}}},
\ldots,\frac{x_n}{M^{\frac1{2\theta}}}\right)=
\pro_{i=1}^n|x_i|^{\alpha}\cdot\det\left[ \Cal K^{Her(\alpha,\theta)}
(x_i,x_j)\right]_{i,j=1}^n
$$ 
where
$$
\Cal K^{Her(\alpha,\theta)}(x,y)=\Cal K^{\left(\frac{\alpha-1}{2},\theta\right)}
(x^2,y^2)+x^{\theta}y\cdot 
\Cal K^{\left(\frac{\alpha+\theta}{2},\theta\right)}(x^2,y^2)
\tag 5.2
$$
and $\Cal K^{(\alpha,\theta)}(x,y)$ is defined in (3.6).
\endproclaim
The proof is straightforward.
\example{Example 5.5} The asymptotic kernel of the
 classic Hermite ensemble in the bulk of spectrum is the {\it sine-kernel}
  $$\frac{\sin\pi(\xi-\eta)}{\pi(\xi-\eta)}.$$
   Let us obtain it from our formulas. For $\theta=1$ we have, see (3.6),
$$
\Cal K^{(\alpha,1)}(x,y)=\int_0^1J_{\alpha+1,{1}}(xt)\cdot
 J_{{\alpha+1},{1}}(yt)t^{\alpha}dt
$$
But 
$$
\Cal K^{Her(0,1)}(x,y)=\Cal K^{\left(-\frac{1}{2},1\right)}(x^2,y^2)+xy\cdot 
\Cal K^{\left(-\frac{3}{2},1\right)}(x^2,y^2)
\tag 5.3
$$
and
$$
J_{\frac 12,1}(x)=\sum_{k=0}^{\infty}\frac{(-x)^k}{k!\Gamma(k+1/2)}=\frac 1{\sqrt\pi}\cos(2\sqrt x),
$$$$
J_{\frac 32,1}(x)=\sum_{k=0}^{\infty}\frac{(-x)^k}{k!\Gamma(k+3/2)}=\frac 1{\sqrt \pi}\frac{\sin(2\sqrt x)}{\sqrt{x}}.
$$
Thus,
$$
\multline
\Cal K^{\left(-\frac 12,1\right)}(x^2,y^2)=\frac 1{\pi}\int_0^1 \cos(2x\sqrt {t})\cos(2y\sqrt {t})\frac{dt}{\sqrt{t}}\\ =\frac 1{2\pi}\left(\frac{\sin2(x-y)}{x-y}+\frac{\sin2(x+y)}{x+y}\right),
\endmultline
$$
$$
\multline
xy\cdot \Cal K^{\left(\frac 12,1\right)}(x^2,y^2)=\frac 1{\pi}\int_0^1{ \sin(2x\sqrt {t})\sin(2y\sqrt {t})}\frac{dt}{\sqrt{t}}\\ =\frac 1{2\pi}\left(\frac{\sin2(x-y)}{x-y}-\frac{\sin2(x+y)}{x+y}\right).
\endmultline
$$
Then (5.3) brings us to the sine-kernel for $\xi=2x/\pi$, $\eta=2y/\pi$.
\endexample
 Similarly to the Laguerre case, the biorthogonal Hermite ensemble also possesses a strange symmetry, cf. Corollary 4.7.
\proclaim{Corollary 5.6} The asymptotic behaviour of the Hermite biorthogonal
ensemble in the bulk of spectrum
 is invariant with respect to the change of parameters (4.5)
 and the transformation $x\mapsto \operatorname{sign}(x)\cdot |x|^{\frac 1\theta}$ of the phase space $(-\infty,+\infty)$.
\endproclaim
\demo{Proof} The claim easily follows from Corollary 4.7 and (5.2). 
If we set, see (4.5),
 $$
\tilde\alpha = \frac{\alpha+1}\theta -1,\quad \tilde \theta= \frac 1\theta
$$ 
then 
$$
\frac{\frac{\alpha-1}{2}+1}\theta -1=\frac{\tilde\alpha-1}2,\quad \frac{\frac{\alpha+\theta}{2}+1}\theta -1=\frac{\tilde\alpha+\tilde\theta}{2}.
$$
These identities show that each summand of (5.2) is invariant under changes from the hypothesis. \qed
\enddemo
It would be very interesting to find some kind of  natural explanation
 for Corollaries 4.7 and 5.6.
\head 6. Appendix\endhead
\subhead Proof of Theorem 3.4\endsubhead
The formula (3.5) for the Christoffel-Darboux kernel $K_N^{Jac}(x,y)$ implies
that
$$
K_N^{Jac}(x,y)=\theta\inte_{0}^1 A_N(xt)B_N((yt)^{\theta})t^{\alpha}dt
\tag 6.1
$$
where
$$
A_N(x)=
\su_{k=1}^N\frac{\left(\frac{k+\alpha}{\theta}\right)_N(-x)^{k-1}}{(k-1)!(N-k)!}
\tag 6.2
$$
$$
B_N(y)=\su_{l=1}^N\frac{
(\theta (l-1)+\alpha+1)_N(-y)^{l-1}}{(l-1)!(N-l)!}
\tag 6.3
$$
Comparing (6.1) to (3.6) we see that it suffices to prove the following
$$
\lim\limits_{N\to\infty}\frac 1{N^{\frac{\alpha+1}{\theta}}}A_N\left(\frac
x{N^{1+\frac1\theta}}\right)=J_{\frac{\alpha+1}{\theta},\frac{1}\theta}(x),
\tag 6.4
$$
$$
\lim\limits_{N\to\infty}\frac 1{N^{\alpha+1}}B_N\left(\frac
y{N^{1+\frac1\theta}}\right)=J_{{\alpha+1},{\theta}}(y)
\tag 6.5
$$
where the convergence is uniform on every compact subset of $\Bbb R_+$. 
Indeed, these relations imply
$$
\lim\limits_{N\to\infty}\frac 1{N^{1+\frac{1}{\theta}}}\left(\frac{x}{N^{1+\frac{1}{\theta}}}\right)^{\alpha}
K_N^{Jac}\left(\frac{x}{N^{1+\frac{1}{\theta}}}
,\frac{y}{N^{1+\frac{1}{\theta}}}\right)=x^{\alpha}\Cal K^{(\alpha,\theta)}(x,y),
$$
which is what we want to prove.

We shall prove only (6.4); the proof of (6.5) is 
quite similar.

Let us split the sum (6.2) into two parts: in the first part the summation
index $k$ runs from $1$ to some $M<N$ which will be chosen later,
  and the second part is the
 remainder. We
shall denote these parts by $A_N'(x)$ and $A_N''(x)$, respectively.
Thus,
$$
A_N(x)=A_N'(x)+A_N''(x).
$$
It will be sufficient to prove the following (uniform) estimates
$$
\frac 1{N^{\frac{\alpha+1}{\theta}}}A_N'\left(\frac
x{N^{1+\frac1\theta}}\right)=\su_{k=0}^{M-1}\frac{(-1)^{k}x^{k}}
{k!\Gamma(\frac{\alpha+1+ k}{\theta} )}+o(1);
\tag 6.6
$$
$$
\frac 1{N^{\frac{\alpha+1}{\theta}}}A_N''\left(\frac
x{N^{1+\frac1\theta}}\right)=o(1)
\tag 6.7
$$
as $N\to \infty$.

To verify (6.6) we shall use Stirling formula
$$
\ln \Gamma(z)=\left(z-\frac 12\right)\ln z - z +\frac 12 \ln(2\pi)+o(z^{-1}).
\tag 6.8
$$
Applying (6.8) to $z=N+\frac{k+\alpha}{\theta}$ and to $z=N-k+1$ we get (we
shall choose $M$ so that it will be $o(N)$, that is why the last term in (6.8)
produces $O(\frac 1N)$ in the next formula)
$$
\multline
\ln\frac{\Gamma\left(N+\frac{k+\alpha}{\theta}\right)}
{\Gamma(N-k+1)N^{\frac{k+\alpha}{\theta}+k-1}}=
\left(N+\frac{k+\alpha}{\theta}-\frac12\right)
\ln\left(N+\frac{k+\alpha}{\theta}\right)\\-\left(N-k+\frac12\right)\ln(N-k+1)
-\left(\frac{k+\alpha}{\theta}+k-1\right)(\ln N+1)+O\left(\frac1N\right).
\endmultline
$$
Using asymptotic expansions
$$
\ln\left(N+\frac{k+\alpha}{\theta}\right)=\ln N+\frac{k+\alpha}{N\theta}+
O\left(\left(\frac{k+\alpha}{N\theta}\right)^2\right),
$$
$$
\ln(N-k+1)=\ln N-\frac{k-1}{N}+O\left(\left(\frac{k-1}{N}\right)^2\right)
$$ 
we arrive at the following estimate
$$
\multline
\ln\frac{\Gamma\left(N+\frac{k+\alpha}{\theta}\right)}
{\Gamma(N-k+1)N^{\frac{k+\alpha}{\theta}+k-1}}=
\frac{k+\alpha}{N\theta}\left(\frac{k+\alpha}{\theta}-\frac12\right)\\+
\frac{k-1}{N}\left(-k+\frac12\right)+O\left(\frac {k^2}N\right)+
O\left(\frac1N\right)
=O\left(\frac {k^2}N\right).
\endmultline
$$
Now we want this expression to converge to $0$ as $N\to \infty$.
Since $k\leq M$, we may set, for example, $M=[N^{\frac13}]$. Then 
$$
O\left(\frac {k^2}N\right)=O(N^{-\frac13}),
$$
and hence we get
$$
\frac{\Gamma\left(N+\frac{k+\alpha}{\theta}\right)}
{\Gamma(N-k+1)N^{\frac{k+\alpha}{\theta}+k-1}}=1+O(N^{-\frac13} )
$$
and
$$
\multline
\frac 1{N^{\frac{\alpha+1}{\theta}}}A_N'\left(\frac
x{N^{1+\frac1\theta}}\right)=\frac 1{N^{\frac{\alpha+1}{\theta}}}\su_{k=1}^{M}
\frac{\left(\frac{k+\alpha}{\theta}\right)_N}{(k-1)!(N-k)!}\left(\frac{-x}
{N^{1+\frac1\theta}}\right)^{k-1}\\ =
\su_{k=1}^{M}
\frac{\Gamma\left(N+\frac{k+\alpha}{\theta}\right)}
{\Gamma(N-k+1)N^{\frac{k+\alpha}{\theta}+k-1}}\frac{(-x)^{k-1}}{(k-1)!\Gamma\left(
\frac{k+\alpha}{\theta}\right)}\\=(1+O(N^{-\frac13}))\su_{k=1}^{M}
\frac{(-x)^{k-1}}{(k-1)!\Gamma\left(
\frac{k+\alpha}{\theta}\right)}=\su_{k=0}^{M-1}\frac{(-1)^{k}x^{k}}
{k!\Gamma(\frac{\alpha+1+ k}{\theta} )}+o(1).
\endmultline
$$
Thus, (6.6) is proved. 

To prove (6.7) we notice that for any $a,b>0$, $b$  is an integer $<N$, 
we have the following simple
estimate
$$
\frac{\Gamma(N+a)}{\Gamma(N-b)N^{a+b}}\leq \frac{(N+a)^{a+b+1}}{N^{a+b}}\leq
N\left( 1+\frac aN \right)^{a+b+1},
$$
which can be obtained by applying applying the identity
$\Gamma(z+1)=z\Gamma(z)$ to the numerator of the left-hand side $[a]+b+1$
times if $N-b>1$, and $[a]+b$ times if $N-b=1$.

Then
$$
\multline
\left|\frac 1{N^{\frac{\alpha+1}{\theta}}}A_N''\left(\frac
x{N^{1+\frac1\theta}}\right)\right|=\left|\frac 1{N^{\frac{\alpha+1}{\theta}}}
\su_{k=M+1}^{N}
\frac{\left(\frac{k+\alpha}{\theta}\right)_N}{(k-1)!(N-k)!}\left(\frac{-x}
{N^{1+\frac1\theta}}\right)^{k-1}\right|\\ \leq\su_{k=M+1}^{N}
\frac{\Gamma\left(N+\frac{k+\alpha}{\theta}\right)}
{\Gamma(N-k+1)N^{\frac{k+\alpha}{\theta}+k-1}}\frac{|x|^{k-1}}{(k-1)!
\Gamma\left(
\frac{k+\alpha}{\theta}\right)}\\ \leq\su_{k=M+1}^{N}
N\left(1+\frac{k+\alpha}{N\theta}\right)^{\frac{k+\alpha}{\theta}+k}
\frac{|x|^{k-1}}{(k-1)!
\Gamma\left(
\frac{k+\alpha}{\theta}\right)}
\endmultline
$$
 Since $1+\frac{N+\alpha}{N\theta}$ is bounded, say, by some constant $c$, and
 $N\leq (M+1)^3\leq k^3$, the last sum does not
 exceed
 $$
\su_{k=M+1}^{N} k^3
c^{\frac{k+\alpha}{\theta}+k}
\frac{|x|^{k-1}}{(k-1)!
\Gamma\left(
\frac{k+\alpha}{\theta}\right)},
 $$
 which is the difference of
 two partial sums $S_N$ and $S_{M}$ of the
 series
 $$
 \su_{k=1}^{\infty} k^3
c^{\frac{k+\alpha}{\theta}+k}
\frac{|x|^{k-1}}{(k-1)!
\Gamma\left(
\frac{k+\alpha}{\theta}\right)}.
$$
This series converges for all $x$ uniformly on every compact set, and,
consequently,
$$
\frac 1{N^{\frac{\alpha+1}{\theta}}}A_N''\left(\frac
x{N^{1+\frac1\theta}}\right)
$$
converges to $0$ uniformly on every compact set, as was to be proved. \qed
\subhead Proof of Theorem 4.5 \endsubhead 
Let us rewrite the formula (4.4) for the kernel
$K_N^{Lag}(x,y)$ in the following form
$$
\aligned
K_N^{Lag}(x,y)=\theta\sum_{k,i=0}^{N-1}
\frac{\Gamma(N)}{\Gamma(\alpha+\theta k+1)k!(N-k-1)!}
\frac{\Gamma\left(N+\frac{i+\alpha+1}{\theta}\right)}
{\Gamma(N)\Gamma\left(\frac{i+\alpha+1}{\theta}\right)i!}&
\\ \times\frac{(-x^{\theta})^k(-y)^i}{\alpha+\theta k
+i+1}\sum_{r=i}^{N-1}&\frac{y^{r-i}}{(r-i)!}
\endaligned
\tag 6.8
$$
Note now, that if we substitute ${y}{N^{-\frac 1\theta}}$ instead of $y$ into
the last sum, it will be close to $1$:
$$
\sum_{r=i}^{N-1}\frac{1}{(r-i)!}\left(\frac y{N^{\frac
1\theta}}\right)^{r-i}=1+O(N^{-\frac
1\theta})
$$
where $O(N^{-\frac1\theta})$ does not depend on $i$.
Thus, we can neglect this sum while computing the limit.

The rest of (6.8) can be written in the form (cf. (6.1)) 
$$
\aligned
\theta\sum_{k,i=0}^{N-1}
\frac{\Gamma(N)}{\Gamma(\alpha+\theta k+1)k!(N-k-1)!}
&\frac{\Gamma\left(N+\frac{i+\alpha+1}{\theta}\right)}
{\Gamma(N)\Gamma\left(\frac{i+\alpha+1}{\theta}\right)i!}
\frac{(-x^{\theta})^k(-y)^i}{\alpha+\theta k
+i+1}\\ &=
\theta \inte_{0}^1C_N((xt)^{\theta})D_N(yt)t^{\alpha}dt
\endaligned
\tag 6.9
$$
where
$$
C_N(x)=\sum_{k=0}^{N-1}
\frac{\Gamma(N)(-x)^k}{\Gamma(\alpha+\theta k+1)k!(N-k-1)!},
$$
$$
D_N(y)=\sum_{i=0}^{N-1}\frac{\Gamma\left(N+\frac{i+\alpha+1}{\theta}\right)(-y)^i}
{\Gamma(N)\Gamma\left(\frac{i+\alpha+1}{\theta}\right)i!}.
$$
Comparing (6.9) to (3.6) we see that it suffices to prove the following 
$$
\lim\limits_{N\to\infty}C_N\left(\frac
x{N^{\frac1\theta}}\right)=J_{{\alpha+1},1}(x),
\tag 6.10
$$
$$
\lim\limits_{N\to\infty}\frac 1{N^{\frac {\alpha+1}\theta}}D_N\left(\frac
y{N^{\frac1\theta}}\right)=J_{\frac{\alpha+1}{\theta},\frac1{\theta}}
(y)
\tag 6.11
$$
as $N\to\infty$, because these relations imply the desired one
$$
\lim\limits_{N\to\infty}\frac 1{N^{\frac{1}{\theta}}}
\left(\frac{x}{N^{\frac{1}{\theta}}}\right)^{\alpha}\exp\left(\frac x{N^{\frac
1\theta}}\right)
K_N^{Lag}\left(\frac{x}{N^{\frac{1}{\theta}}}
,\frac{y}{N^{\frac{1}{\theta}}}\right)=x^{\alpha}\Cal K^{(\alpha,\theta)}(y,x).
$$
(The interchange $x\leftrightarrow y$ in the last expression does not change the
determinants of the type $\det [\Cal K^{(\alpha,\theta)}(x_i,x_j)]$.)

The proofs of (6.10) and (6.11) are very similar to the proof
of (6.4) which we carried out above, and we shall not give them here. \qed

\bigskip
 Department of Mathematics, The University of
Pennsylvania, Philadelphia, PA 19104-6395, U.S.A.  E-mail address:
{\tt borodine\@math.upenn.edu} 
\enddocument